\documentclass[11pt]{article}
\usepackage{amssymb}
\usepackage{amsfonts}
\usepackage{amsthm}
\usepackage{graphicx}
\usepackage[utf8]{inputenc}
\usepackage{xcolor}
\usepackage[english]{babel}
\usepackage{amsmath}
\usepackage{eucal}
\usepackage{enumerate}
\usepackage{stmaryrd}
\usepackage{amsfonts,epsfig,epstopdf,titling,url,array}
\usepackage{verbatim} 
\usepackage[multiple]{footmisc}

\numberwithin{equation}{section}

\newtheorem{thm}{Theorem}
\newtheorem{prop}{Proposition}
\newtheorem{lemma}{Lemma}

\newcommand{\beq}{\begin{equation}}
\newcommand{\RN}[1]{%
  \textup{\uppercase\expandafter{\romannumeral#1}}%
}
\newcommand{\eeq}{\end{equation}}

\newcommand{\beqs}{\begin{equation*}}
\newcommand{\eeqs}{\end{equation*}}

\usepackage[
top    = 2.75cm,
bottom = 2.50cm,
left   = 2.49cm,
right  = 2.549cm]{geometry}

\title{The effects of simple density-dependent prey diffusion and refuge in a predator-prey system}

\author{
Leoncio Rodriguez Q., Jia Zhao, Luis F. Gordillo\\ Department of Mathematics and Statistics \\ Utah State University \\
Logan, UT 84322.
}

\begin{document}
\bibliographystyle{plain}
\large{

\begin{titlepage}
\LARGE{

\vskip 1truein
\begin{center}
The effects of simple density-dependent prey diffusion and refuge in a predator-prey system
\end{center}
}
\large{
\vskip 1truein
\begin{center} 
Leoncio Rodriguez Q.\footnote{ Corresponding author, e-mail address: leoncio.quinones@usu.edu}, Jia Zhao\footnote{e-mail address: jia.zhao@usu.edu}, Luis F. Gordillo.\footnote{e-mail address:  luis.gordillo@usu.edu} \\
Department of Mathematics and Statistics\\
Utah State University \\
Logan, UT 84322\\
\end{center}

\vskip 1truein
\begin{abstract}
We study a spatial (two-dimensional) Rosenzweig-MacArthur model under the following assumptions: $(1)$ prey spread follows a nonlinear diffusion rule, $(2)$ preys have a refuge zone (sometimes called ``protection zone") where predators cannot enter, (3) predators move following linear diffusion. We present a bifurcation analysis for the system that shows the existence of positive solutions at the steady state. We complement the theoretical results with numerical computations and compare our results with those obtained in the case of having linear diffusion for the prey movement. Our results show that both models, with linear and nonlinear diffusion for the prey, have the same bifurcation point and the positive solution curves are virtually the same in a neighborhood of this point, but they get drastically different as the bifurcation parameter approaches to zero. 

\end{abstract}
}
\end{titlepage}

\maketitle

\section{Introduction}
It seems completely natural to expect the effects of density-dependent dispersal on spatially distributed predator-prey systems \cite{Okubo}; for instance, reduced amounts of resources at fixed spatial locations due to a high level of aggregation might drive the dispersal for searching and acquiring new resources for survival at locations with less competition. Other adaptive responses could also be observed in some prey and predator populations, like keeping away from crowds to be less visible to predators and avoiding encounters with conspecific individuals in active searching of prey to decrease interference \cite{Hassell}. 

A theoretical framework for density-dependence dispersal includes diffusivity as a function of the population density in a reaction-diffusion equation, which could also contain additional nonlinear terms regarding other relevant aspects of the system. Here we are interested in studying the effects of nonlinear diffusion by the prey under two specific circumstances: (1) there is predator saturation on prey consumption (we use a Holling type II function) and (2) the prey habitat contains a refuge zone where their predation is not possible and can be thought as a mechanism for conservation \cite{sih}. For the nonlinearity in the diffusion, we assume the simple form $\nabla\cdot u\nabla u$ ($u$ represents the prey population), which is a particular case of a more general model discussed in \cite{Okubo}, see also \cite{Gurtin}.

Although there are variants of the model presented here that have been extensively studied in recent years, see for instance \cite{He2017} and \cite{Du}, we have not found in the literature results that directly compare the effects of density-dependence dispersal with those from linear diffusion, under the conditions (1) and (2) mentioned above. To understand how the differences in the dynamics depend on the model parameters might become relevant when attempting the modeling in real scenarios, as could be in the case of pest suppression efforts through biological control.

We start by showing the existence of nontrivial solutions in the steady state via bifurcation analysis and then we compare numerically the effects of the nonlinearity in the diffusion with its linear counterpart. There are studies involving the simultaneous effects of nonlinearities in the reaction part (in particular, the predator saturation) and refuge, see for instance \cite{wang-li} and \cite{zhangrongzhang}, but the introduction of nonlinear diffusion requires the development of alternative theoretical tools. 

Our particular model of interest is defined over a bounded domain $\Omega\subset\mathbb{R}^2$, which is the representation of a closed environment where predators and preys live. We consider an additional domain, the ``refuge zone", $\Omega_{0}\subset\Omega$, where  predators cannot enter. We assume that $\Omega$  and $\Omega_0$ have sufficiently smooth boundaries, that $\overline{\Omega}_{0}\subset\Omega$, and define $\Omega_{1}=\Omega\setminus\overline{\Omega}_{0}$. Let us consider the following system of parabolic equations for the prey and predator populations, denoted by $u$ and $v$ respectively, 
\begin{equation}\label{system0}
\begin{aligned}
    \partial_{t}u &= D_u\nabla\cdot u\nabla u + ru\left(1 - \frac{u}{\lambda}\right) - \frac{b(x)uv}{1+mu}\qquad\text{in}\quad \Omega, \\ 
    \partial_{t}v &= D_v\Delta v - \mu v + \frac{cuv}{1+mu} \qquad\qquad \text{in}\quad \Omega_{1},\\
    v&\equiv 0\qquad \text{in}\quad \Omega\setminus\Omega_1,
\end{aligned}
\end{equation}
with boundary and initial conditions given by
\begin{equation}
\begin{aligned}
    \partial_{n}u &= 0 \quad\text{on}\quad \partial\Omega, \\
    \partial_{n}v &= 0 \quad\text{on}\quad \partial\Omega_{1},\\ 
    u(x,0) &= u_{0}(x)\geq 0 \quad\text{for}\quad x \in \Omega, \\
    v(x,0) &= v_{0}(x)\geq 0 \quad\text{for}\quad x \in \Omega_{1}. 
\end{aligned}
\end{equation}
The reaction part of this system is the well-known Rosenzweig-MacArthur model, \cite{RosenzweigMac, Kot}, where the parameters are positive and the function $b(x)$, which determines the efficiency of predator attacks, is defined by
\begin{equation}
b(x) = 
    \begin{cases}
    b > 0& \qquad\text{if } x \in \Omega_{1},\\
    0& \qquad \text{if } x \in \overline{\Omega}_{0},
    \end{cases}
\end{equation}
thus characterizing the refuge zone $\Omega_{0}$. By imposing a non-flux boundary condition on $\partial\Omega_{1}$, we restrict predators to the exterior of the refuge zone. In contrast, preys can move freely over the whole domain $\Omega$.
For the bifurcation analysis in the next Section it is convenient to consider the model in dimensionless form. After a suitable re-scaling the two equations in (1.1) can be rewritten as  
\begin{equation}\label{system1}
\begin{aligned}
    \partial_{t}u &= \nabla\cdot u\nabla u +\lambda u - u^{2} - \frac{b(x)uv}{1+mu}\qquad\text{in}\quad \Omega, \\ 
    \partial_{t}v &= d \Delta v - \mu v + \frac{cuv}{1+mu}\qquad\text{in}\quad \Omega_{1},\\ 
\end{aligned}
\end{equation}
where $d=D_v/D_u$. Please notice that, although we are using the initial notation, the variables and parameters have now different interpretation.
%
First we focus on showing the existence of positive steady-state solutions for the homogeneous system
\begin{equation}\label{system2}
\begin{aligned}
    \nabla\cdot u\nabla u +\lambda u - u^{2} - \frac{b(x)uv}{1+mu} &= 0 \qquad\text{in}\quad \Omega \\ 
    \Delta v - \mu v + \frac{cuv}{1+mu} &= 0\qquad\text{in}\quad \Omega_{1}\\ 
    \partial_{n}u &= 0\qquad \text{on}\quad \partial{\Omega}, \\ \partial_{n}v &= 0 \qquad\text{on}\quad \partial\Omega_{1},\\ 
    \end{aligned}
\end{equation}
where the parameters in the predators equation have been redefined accordingly.

\section{Bifurcation Analysis}\label{sec2}
In this Section we show the emergence of positive solutions for the problem (\ref{system2}) and its counterpart that has the laplacian $\Delta u$ for the prey equation replacing the term $\nabla\cdot u\nabla u$. This is achieved by using the Crandall-Rabinowitz theorem on bifurcations from simple eigenvalues. In what follows we assume that $u$ is bounded away from zero.

\subsection{The nonlinear diffusion case}
The first step is to establish the nature of the non-negative solutions, which is done in Proposition 1. For its proof we first require a Lemma that adapts a maximum principle in \cite{lou} to the case of the nonlinear diffusion considered here.
\begin{lemma}\label{lema1}
Suppose $g\in C(\overline{\Omega}\times\mathbb{R})$ and $u\in C^{2}(\Omega)\cap C^{1}(\overline{\Omega})$, $u\geq 0$ in $\Omega$, where $\Omega$ is a bounded domain in $\mathbb{R}^{N}$ with smooth boundary. 
\begin{enumerate}[i.]
   \item If $\nabla\cdot u\nabla u + g(x,u(x))\geq 0$ in $\Omega$, $\partial_{n}u\leq 0$ on $\partial\Omega$ and $u(x_{0})=\max_{\overline{\Omega}}u(x)$, then $g(x_{0},u(x_{0}))\geq 0$
   \item If $\nabla\cdot u\nabla u + g(x,u(x))\leq 0$ in $\Omega$, $\partial_{n}u\geq0$ on $\partial\Omega$ and $u(x_{0})=\min_{\overline{\Omega}}u(x)$, then $g(x_{0},u(x_{0}))\leq 0$
\end{enumerate}
\end{lemma}

\begin{proof}
Part (i). Notice that by continuity of $u$ and compactness of $\overline{\Omega}$, there exists $x_{0}\in\overline{\Omega}$ such that $u(x_{0})=\max_{\overline{\Omega}}u(x)$. If $x_{0}\in\Omega$, then we must have $\Delta u(x_{0})\leq 0$, and $\nabla u(x_{0})=\Vec{0}$. Since $u\geq 0$ in $\Omega$, we have $$(\nabla\cdot u\nabla u)(x_{0})+g(x_{0},u(x_{0}))=u(x_{0})\Delta u(x_{0})+g(x_{0},u(x_{0}))\geq 0.$$
From this, we obtain $0\geq u(x_{0})\Delta u(x_{0})\geq - g(x_{0},u(x_{0}))$, hence $g(x_{0},u(x_{0}))\geq 0$. Now suppose that $x_{0}\in\partial\Omega$ and $g(x_{0},u(x_{0}))<0$. By the continuity of $g$ and $u$, there exists a ball $\textbf{B}\subset\overline{\Omega}$ such that $\partial\Omega\cap\partial \textbf{B} = \{x_{0}\}$. By the hypothesis we have
\begin{equation}\label{equ1}
    u\Delta u + |\nabla u(x)|^{2}\geq -g(x,u(x))>0.
\end{equation}
Let us write $u(x_{0})=\max_{\overline{\Omega}}u(x)\equiv M$ and $v=v(x)\equiv M$. Notice that the term $u\Delta u + |\nabla u|^{2}$ has the form $a_{ij}(x,u,\nabla u)\partial_{x_{ij}}^{2}u+\mathcal{B}(x,u,Du)$ with $a_{ij}=0$ if $i\neq j$, $a_{ii}(x,u,\nabla u) = u$ and $\mathcal{B}(x,u,\nabla u) = |\nabla u|^{2}$. Furthermore, since $u\in C^{2}(\Omega)\cap C^{1}(\overline{\Omega})$, the matrix $[a_{ij}]$ is continuous and continuously differentiable with respect to its second and third arguments in the set $\Omega\times\mathbb{R}\times\mathbb{R}^{N}$, also $\mathcal{B}(x,z,\zeta)=|\zeta|^{2}$, $\zeta\in\mathbb{R}^{n}$ is continuously differentiable with respect to $\zeta$ in $\mathbb{R}^{N}$. The inequality (\ref{equ1}) implies that $u$ is an elliptic solution in the sense described in \cite{pucciserrin} (section 2.2), as well as $v\equiv M$, since $v\Delta v + |\nabla v|^{2}\leq 0$. Let $K$ be a compact subset of  $\Omega\times\mathbb{R}\times\mathbb{R}^{N}$ and choose $\Tilde{z}>z$, the inequality $\mathcal{B}(x,\Tilde{z},\zeta)-\mathcal{B}(x,z,\zeta)\geq-\kappa(\Tilde{z}-z)$ for some $\kappa>0$, follows from the fact that $\mathcal{B}(x,\Tilde{z},\zeta)-\mathcal{B}(x,z,\zeta)=0>-(\Tilde{z}-z)$ with $\kappa=1$, hence $\mathcal{B}$ is lower Lipschitz continuous in the variable $z$ in $K$. Notice that $u<v$ in $\Omega$ and $u=v$ exactly at $x_{0}\in\partial\Omega$, then by Theorem $2.7.1$ in \cite{pucciserrin} we have $\partial_{n}u(x_{0})>\partial_{n}v(x_{0})\equiv 0$, contradicting the boundary condition $\partial_{n}u(x_{0})\leq 0$. Therefore, we must have $g(x_{0},u(x_{0}))\geq 0$ as needed. Part (ii) of the Lemma is proved by a similar argument but reversing the inequalities.
\end{proof}

\begin{prop}\label{propo1}
Let $u,v\in C^{2}(\Omega)\cap C^{1}(\overline{\Omega})$, the non-negative trivial and semi-trivial solutions $(u,v)$ of the system (\ref{system2}), is either $(0,0)$ or $(\lambda,0)$ respectively. 
\end{prop}
\begin{proof}
Suppose $v\equiv 0$, then we want to show that either $u\equiv 0$ or $u\equiv\lambda$. Suppose $u\not\equiv 0$, then $u>0$ in some bounded subset $A$ of $\Omega$ of positive measure. By continuity of $u$, there exists $x_{0}\in\overline{\Omega}$ such that $0<u(x_{0})=\max_{\overline{\Omega}}u(x)$. We have that $\partial_{n}u \leq 0$ on  
$\partial \Omega$  and
$u\Delta u + |\nabla u|^{2}+u(\lambda-u)\geq0$ in $\Omega$, then by Lemma \ref{lema1}, $u(x_{0})(\lambda-u(x_{0}))\geq 0$, which implies that $\lambda\geq u(x_{0})\geq u$ in $\Omega$. On the other hand, since $u$ is continuous, there exists $x^{\prime}\in \overline{\Omega}$ such that $u(x^{\prime})=\min_{\overline{\Omega}}u(x)$. We also have $u\Delta u + |\nabla u|^{2}+u(\lambda-u)\leq0$ in $\Omega$ and $\partial_{n}u\geq 0$ on  
$\partial \Omega$, then by Lemma \ref{lema1} we conclude $u(x^{\prime})(\lambda-u(x^{\prime}))\leq 0$. If $u(x^{\prime})=0$, then the last inequality would imply  $\lambda \leq 0$ which is not possible since we assume $\lambda>0$, thus we must have $u(x^{\prime})>0$ and $\lambda\leq u (x^{\prime})\leq u $ in $\Omega$. Finally, we have $\lambda\leq u \leq \lambda$ in $\Omega$, and therefore $u\equiv\lambda$ in $\Omega$ as needed. 
\end{proof}
Proposition \ref{propo1} is the first step to study the presence of a bifurcation along the semi-trivial solutions given by the curve $\Gamma_{u}=\{(\mu,u,v)=(\mu,\lambda,0):\mu>0\}$, see \cite{He2017} and \cite{Du}. We can now proceed to investigate positive solutions to the system (\ref{system2}).
Although the analysis below follows the mathematical framework presented in \cite{He2017}, \cite{Du}, and  Section 3.4.2 of \cite{cantrell}, we have tried to provide a more detailed account  with the hope of facilitating the reading.
\\
\\
\noindent
Let us start defining the spaces
\begin{align*}
 X_{\Omega} &=\{u\in W^{2,p}(\Omega):\partial_{n}u=0 \text{ on } \partial\Omega\}, \  Y_{\Omega}=L^{p}(\Omega),\\
 X_{\Omega_{1}} &= \{u\in W^{2,p}(\Omega_{1}):\partial_{n}u=0 \text{ on } \partial\Omega_{1}\},\ Y_{\Omega_{1}}=L^{p}(\Omega_{1}). 
\end{align*}
where we assume $p>2$ so that Theorem 3.3 (together with Remark 3.4.1) in \cite{shiwang} can be applied. We are interested in analyzing how nontrivial solutions, i.e. $v>0$, bifurcate in a neighborhood of $(\mu,\lambda,0)$. Consider the function $w=\lambda-u$ and define the operator $F:\mathbb{R}\times X_{\Omega}\times X_{\Omega_{1}}\rightarrow Y_{\Omega}\times Y_{\Omega_{1}}$ by 
\begin{equation}\label{operator1}
F(\mu,w,v)=\begin{pmatrix}
-\nabla\cdot w\nabla w + \lambda\Delta w - \lambda w + w^{2} + \frac{b(x)(\lambda-w)v}{1+m(\lambda-w)}\\ \\
\Delta v - \mu v + \frac{c(\lambda - w)v}{1+m(\lambda-w)}
\end{pmatrix}^{T}.   
\end{equation}
Since we aim to use Theorem 1.7 in \cite{crandall}, we linearize the operator in (\ref{operator1}) by computing 
\[F_{(w,v)}(\mu,w,v)[\alpha,\beta]= \frac{d}{d\epsilon}F(\mu,w+\epsilon\alpha,v+\epsilon\beta) |_{\epsilon=0}.\]
This gives
\begin{equation}\label{operator2}
\resizebox{.9\hsize}{!}{$
F_{(w,v)}(\mu,w,v)[\alpha,\beta]=\begin{pmatrix}
-\nabla\cdot \alpha\nabla w -\nabla\cdot w\nabla\alpha + \lambda\Delta \alpha - \lambda\alpha + 2w\alpha + \frac{b(x)\beta(\lambda-w)}{1+m(\lambda-w)}-\frac{b(x)\alpha v}{(1+m(\lambda-w))^{2}} \\ \\
\Delta \beta - \mu \beta + \frac{c(\lambda - w) \beta}{1+m(\lambda - w)}-\frac{cv\alpha}{(1+m(\lambda-w))^{2}}
\end{pmatrix}^{T}.$}
\end{equation}
Around $(\lambda,0)$, i.e $\lambda = u$ and $v=0$, we have that (\ref{operator2}) becomes

\begin{equation}\label{operator3}
    F_{(w,v)}(\mu,0,0)[\alpha,\beta]=\begin{pmatrix}
    \lambda\Delta\alpha - \lambda\alpha + \frac{b(x)\lambda\beta}{1+m\lambda}\\ \\
    \Delta\beta - \mu\beta + \frac{c\lambda\beta}{1+m\lambda}
    \end{pmatrix}^{T}.
\end{equation}
Since we are interested in nontrivial solutions to (\ref{system2}) we look for the values of $\mu$ for which $F_{(w,v)}(\mu,0,0)[\alpha,\beta]=(0,0)$ has no trivial solutions and $\dim(\text{\textit{ker}}\left( F_{(w,v)}(\mu,0,0)\right))=1$. Notice that $(\alpha,\beta)=(0,0)$ is always a solution to $F_{(w,v)}(\mu,0,0)[\alpha,\beta]=(0,0)$. First We look if $(\alpha,0)$ with $\alpha\neq 0$ is also a solution. If that is the case, then $\lambda\Delta\alpha -\lambda\alpha = 0$ in $\Omega$ and $\partial_{n}\alpha = 0$ on $\partial\Omega$. Then, since $\lambda>0$,  the weak formulation of the associated partial differential equation gives $-\int_{\Omega}\nabla\alpha\cdot\nabla\psi dx = \int_{\Omega}\alpha\psi dx$. For any $\psi\in X_{\Omega}$, in particular for $\psi=\alpha$, we get $-\int_{\Omega}|\nabla\alpha|^{2}dx=\int_{\Omega}\alpha^{2}dx$, which holds only if $\alpha=0$. Thus, we look at solutions of the form $(0,\beta)$ with $\beta\neq0$.

We recall that the Neumann eigenvalues for the Laplacian can be characterized by the $\min$-$\max$ formula, see \cite{courant} or \cite{AHen06} for instance,
\begin{equation}\label{minmaxform}
    \mu^{N}_{k}(\Omega)=\min_{S_{k}\subset H^{1}(\Omega)}\max_{\phi\in S_{k},\phi \neq 0}\frac{\int_{\Omega}|\nabla \phi|^{2}dx}{\int_{\Omega}\phi^{2}dx},
\end{equation}
where $S_{k}$ are subspaces of dimension $k$ of the Sobolev space $H^{1}(\Omega)$ and the minimum is achieved by choosing $S_{k}$ to be the subspace spanned by the first $k$ eigenfunctions $\phi_{1},\phi_{2},\dots,\phi_{k}$.
Notice that $\mu^{N}_{1}(\Omega)=0$ is a consequence of (\ref{minmaxform}), corresponding to the constant eigenfunction $\phi_{1}$ (constant at least on a connected component of $\Omega$). This can be justified as follows. First, notice that $\mu^{N}_{k}(\Omega)\geq 0$ and that zero is achieved whenever $\max_{\phi\in S_{k},\phi \neq 0}\left( \int_{\Omega}|\nabla \phi|^{2}dx / \int_{\Omega}\phi^{2}dx\right)=0$, which is possible only if $\int_{\Omega}|\nabla \phi|^{2}dx / \int_{\Omega}\phi^{2}dx=0$, and hence we must have $\int_{\Omega}|\nabla \phi|^{2}dx = 0$, which implies that $|\nabla \phi| = 0$ a.e on $\Omega$, if $\Omega$ is connected then we have $\phi$ is constant a.e on $\Omega$.

Now, consider the boundary value problem determined by the second component of (\ref{operator3}),
\begin{equation}\label{equ2-6}
\begin{aligned}
    -\Delta\beta &= \left(-\mu+\frac{c\lambda}{1+m\lambda}\right)\beta \ \ \text{in} \ \ \Omega_{1}\\
    \partial_{n}\beta &= 0 \ \ \text{on    } \ \ \partial\Omega_{1}
\end{aligned}
\end{equation}
From (\ref{minmaxform}) we have that  only when $\mu_{\lambda} = c\lambda/(1+m\lambda)$, $\beta$ does not change sign on $\Omega_{1}$ (since this corresponds to the zero eigenvalue). More precisely, $\beta\geq 0$ implies that $\beta$ is a positive constant. Therefore, $\mu_{\lambda}$ is the unique bifurcation point along $\Gamma_{u}$ from which positive solutions of (\ref{system2}) emerge.

The argument above also shows that $\text{\textit{ker}}\left( F_{(w,v)}(\mu_{\lambda},0,0)\right)=span \{(\alpha_{\mu_{\lambda}},1)\},$ where $\alpha_{\mu_{\lambda}}$ solves the boundary value problem
\begin{align}
    \Delta\alpha -\alpha + \frac{b(x)}{1+m\lambda} &= 0 \ \text{in} \ \ \Omega,\\
    \partial_{n}\alpha &= 0  \ \text{in} \ \ \partial\Omega.
\end{align}
Notice that by choosing $\beta=1$, with $\lambda = \mu_{\lambda}/(c-m\mu_{\lambda})$, we get $\alpha_{\mu_{\lambda}}=(-\Delta + I)^{-1}\left[b(x)/(1+m \lambda)\right]$. On the other hand, if we consider the non-homogeneous problem $F_{(w,v)}(\mu_{\lambda},0,0)[\alpha,\beta]=\left[f(x),g(x)\right]$ then the $\text{\textit{Range}}(F_{(w,v)}(\mu_{\lambda},0,0))=\left\{(f,g)\in Y_{\Omega}\times Y_{\Omega_{1}}|\int_{\Omega_{1}}g(x)dx = 0\right\}$. To see this, from the weak formulation we must have
\begin{align}
-\int_{\Omega}\nabla\alpha\cdot\nabla u dx - \int_{\Omega}\alpha u dx + \int_{\Omega_{1}}b\frac{1}{1+m\lambda}u dx &= \int_{\Omega}\frac{f}{\lambda}u dx \\
-\int_{\Omega_{1}}\nabla\beta\cdot\nabla v dx - \mu_{\lambda}\int_{\Omega_{1}}\beta v dx +\int_{\Omega_{1}}\frac{c\lambda}{1+m\lambda}\beta v dx &= \int_{\Omega_{1}}gvdx,
\end{align}
for any $(u,v)\in X_{\Omega}\times X_{\Omega_{1}}$. In particular, if we choose $(u,v)=(\alpha_{\mu_{\lambda}},1)$, we get the condition $\int_{\Omega_{1}}g(x)dx=0$. Therefore
\begin{align}
    \dim(\text{\textit{ker}}\left( F_{(w,v)}(\mu_{\lambda},0,0)\right))=\text{codim}(\text{\textit{Range}}\left( F_{(w,v)}(\mu_{\lambda},0,0)\right))=1.
\end{align}
Notice also that $F_{\mu}(\mu,w,v)=(0,-v)$ and $F_{\mu(w,v)}(\mu,w,v)[\alpha,\beta]=(0,-\beta)$. Therefore,
\begin{align}\label{derivativemu}
    F_{\mu(w,v)}(\mu_{\lambda},0,0)[\alpha_{\mu_{\lambda}},1]=(0,-1).
\end{align}
In particular,  $(0,-1)\not\in\text{\textit{Range}}\left( F_{(w,v)}(\mu_{\lambda},0,0)\right)$. By the classical result on bifurcations from simple eigenvalues of Crandall and Rabinowitz \cite{crandall}, we conclude that the positive solutions of the system (\ref{system2}) form a smooth curve given by 
\begin{align}\label{branchsolution}
\{(\mu,u,v)=(\mu_{\lambda}(s),\lambda-s\alpha_{\mu_{\lambda}}(x)+o(|s|),s+o(|s|)):s\in (0,a)\}
 \end{align}
for some $a>0$, bifurcating from $\Gamma_{u}$ at $(\mu_{\lambda},\lambda,0)$ and such that $\mu_{\lambda}(0)=c\lambda/(1+m\lambda)$.

Along the branch $(\mu(s),w(s),v(s))$ given by (\ref{branchsolution}) the operator defined in (\ref{operator1}) depends on the variable $s$. Thus we compute $F_{ss}(\mu(s),w(s),v(s))$, which is given by

\begin{equation}\label{operator4}
\begin{pmatrix}
-w_{ss}\Delta w - 2w_{s}\Delta w_{s}-w\Delta w_{ss}-2\nabla w\cdot\nabla w_{ss}-2|\nabla w_{s}|^{2}+\lambda\Delta w_{ss}-\lambda w_{ss}+2w^{2}_{s}+2ww_{ss}-\\
\frac{2mb(x)w^{2}_{s}v(s)}{(1+m(\lambda-w(s)))^{3}}-\frac{b(x)w_{ss}v(s)}{(1+m(\lambda-w(s)))^{2}}-\frac{2b(x)w_{s}v_{s}}{(1+m(\lambda-w(s)))^{2}}+\frac{b(x)(\lambda-w(s))v_{ss}}{1+m(\lambda-w(s))}\\
\\
\Delta v_{ss}-\mu^{\prime\prime}(s)v-2\mu^{\prime}(s)v_{s}-\mu v_{ss}-\frac{2mcw^{2}_{s}v(s)}{(1+m(\lambda-w(s)))^{3}}-\frac{cw_{ss}v(s)}{(1+m(\lambda-w(s)))^{2}}-\frac{2cw_{s}v_{s}}{(1+m(\lambda - w(s)))^{2}}+\\
\frac{c(\lambda-w(s))v_{ss}}{1+m(\lambda-w(s))}
\end{pmatrix}^{T}.   
\end{equation}
Using subscripts to denote the first and second derivative of $v$ and $w$, respectively, and using the fact that $w(s)=\lambda-u(s)=\lambda-(\lambda - s\alpha_{\mu_{\lambda}}(x)+o(|s|))$, $v(s) = s + o(|s|)$, we have that at $s=0$, $v(0)=0$, $v_{s}(0)=1$, $v_{ss}(0)=0$, $w(0)=0$, $w_{s}(0)=\alpha_{\mu_{\lambda}}$, $w_{ss}(0)=0$. Therefore, with $\beta = 1$, the expression  (\ref{operator4}) becomes
\begin{equation}\label{operator5}
F_{ss}(\mu_{\lambda}(0),0,0)=\begin{pmatrix}
-2\nabla\cdot(\alpha_{\mu_{\lambda}}\nabla\alpha_{\mu_{\lambda}}) + 2\alpha_{\mu_{\lambda}}^{2} - \frac{2b(x)}{(1+m\lambda)^{2}}\alpha_{\mu_{\lambda}}\beta \\ \\
-2\mu^{\prime}_{\lambda}(0) - \frac{2c}{(1+m\lambda)^{2}}\alpha_{\mu_{\lambda}}\beta
\end{pmatrix}^{T},   
\end{equation}
where $\mu^{\prime}_{\lambda}$ represents the derivative of $\mu_{\lambda}(s)$ with respect to $s$. 

Furthermore, notice that $F_{(w,v)(w,v)}(\mu_{\lambda},w,v)[\alpha_{\mu_{\lambda}},\beta]^{2}$ is equal to
\begin{equation}\label{operator6}
\begin{pmatrix}
-2\nabla\cdot(\alpha_{\mu_{\lambda}}\nabla\alpha_{\mu_{\lambda}}) + 2\alpha_{\mu_{\lambda}}^{2} - \frac{2b(x)}{(1+m(\lambda-w))^{2}}\alpha_{\mu_{\lambda}}\beta  -\frac{2m b(x)}{(1+m(\lambda-w))^{3}}\alpha^{2}_{\mu_{\lambda}}v \\ \\
\frac{-2c}{(1+m(\lambda - w))^{2}}\alpha_{\mu_{\lambda}}\beta  -\frac{2m c}{(1+m(\lambda-w))^{3}}\alpha^{2}_{\mu_{\lambda}}v
\end{pmatrix}^{T},   
\end{equation}
and at $(\mu_{\lambda},0,0)$,
\begin{equation}\label{operator7}
    F_{(w,v)(w,v)}(\mu_{\lambda},0,0)[\alpha_{\mu_{\lambda}},\beta]^{2} =
    \begin{pmatrix}
    -2\nabla\cdot(\alpha_{\mu_{\lambda}}\nabla\alpha_{\mu_{\lambda}}) + 2\alpha_{\mu_{\lambda}}^{2} - \frac{2b(x)}{(1+m\lambda)^{2}}\alpha_{\mu_{\lambda}}\beta\\ \\
    \frac{-2c}{(1+m\lambda)^{2}}\alpha_{\mu_{\lambda}}\beta
    \end{pmatrix}^{T}.
\end{equation} 
As long as $F_{(w,v)(w,v)}(\mu_{\lambda},0,0)\not\in Range(F_{(w,v)}(\mu_{\lambda},0,0))$ we use formula I.6.3 in \cite{Hansjorg} to compute an explicit expression for $\mu^{\prime}_{\lambda}(0)$. Consider the projection acting on $F_{(w,v)(w,v)}(\mu_{\lambda},0,0)[\alpha_{\mu_{\lambda}},\beta]^{2}$, and $F_{\mu(w,v)}(\mu_{\lambda},0,0)$ respectively and defined by

\begin{align}
    \langle F_{(w,v)(w,v)}(\mu_{\lambda},0,0)[\alpha_{\mu_{\lambda}},\beta]^{2}, \pi_{1}\rangle &= \int_{\Omega_{1}}\frac{-2c}{(1+m\lambda)^{2}}\alpha_{\mu_{\lambda}}\beta dx
    \\
    \langle F_{\mu(w,v)}(\mu_{\lambda},0,0)[\alpha_{\mu_{\lambda}},\beta],\pi_{1} \rangle &=  \int_{\Omega_{1}}(-\beta) dx.
\end{align}{}
Then, 
\begin{align}\label{derivativeeval}
    \mu^{\prime}_{\lambda}(0) &= -\frac{1}{2}\frac{\langle F_{(w,v)(w,v)}(\mu_{\lambda},0,0)[\alpha_{\mu_{\lambda}},\beta]^{2}, \pi_{1}\rangle}
{\langle F_{\mu(w,v)}(\mu_{\lambda},0,0)[\alpha_{\mu_{\lambda}},\beta],\pi_{1} \rangle}\\
&=-\frac{c}{|\Omega_{1}|(1+m\lambda)^{2}}\int_{\Omega_{1}}\alpha_{\mu_{\lambda}}dx < 0.
\end{align}
Notice that the condition $F_{(w,v)(w,v)}(\mu_{\lambda},0,0)\not\in Range(F_{(w,v)}(\mu_{\lambda},0,0))$ guarantees that the integral over $\Omega_{1}$ of the second component of ($\ref{operator7}$) does not vanish, so $\mu^{\prime}_{\lambda}(0)$ will not be zero.

To determine the values of $\mu>0$ for which the system (\ref{system2}) has either a unique positive solution, at least one positive solution, or no positive solutions, we use a unilateral global bifurcation result for Fredholm operators due to Shi and Wang (see Theorem 4.4 of \cite{shiwang}), which is based on a result due to Lopez-Gomez (Theorem 6.4.3 page 188 in \cite{lopezgomez}). This updates Rabinowitz's Theorem 1.27, presented originally in \cite{rabinowitz}. Consider the operator $F_{0}:\mathbb{R}\times X_{\Omega}\times X_{\Omega_{1}}\rightarrow Y_{\Omega}\times Y_{\Omega_{1}}$, associated to (\ref{system2}) and defined by
\begin{equation}
    F_{0}(\mu,u,v)=
    \begin{pmatrix}{}
    u-\lambda\\ \\
    v
    \end{pmatrix}^{T}-
    \begin{pmatrix}{}
    (-N + I)^{-1}(u-\lambda +\lambda u - u^{2} - b(x)uv/(1+mu) )\\
    \\
    (-\Delta + I)^{-1}_{\Omega_{1}}(v - \mu v + cuv/(1+mu))
    \end{pmatrix}^{T},
\end{equation}
where $N=N(u)=\nabla\cdot u\nabla u$. Then $F_{0}(\mu,u,v) = \Vec{0}$ is equivalent to the system (\ref{system2}). Let $\mathcal{S}\subset\mathbb{R}\times X_{\Omega}\times X_{\Omega_{1}}$ be the set of non-negative solutions of (\ref{system2}) determined by (\ref{branchsolution}), and let  $\mathcal{C}$ be the connected component of the set $\mathcal{S}\cup\{(\mu_{\lambda},\lambda,0)\}$ emanating from  $(\mu_{\lambda},\lambda,0)$ such that $\mathcal{C}\subset\{(\mu,u,v)\in\mathbb{R}\times X_{\Omega}\times X_{\Omega_{1}}\setminus\{(\mu_{\lambda},\lambda,0)\}:F_{0}(\mu,u,v)=\Vec{0}\}$. 

From the above computations, all conditions of Theorem 4.3 in \cite{shiwang} are satisfied and it can be seen from equations (\ref{operator2}) and (\ref{derivativemu}) that $F_{(u,v)}(\mu,\lambda,0)$ is continuously differentiable in $\mu$ for $(\mu,\lambda,0)$ in $\mathbb{R}\times X_{\Omega}\times X_{\Omega_1}$. The $C^{1}$ condition of the corresponding norm function $(u,v)\mapsto \|(u,v)\|_{X_{\Omega}\times X_{\Omega_{1}}}$ is guaranteed in \cite{restrepo} (see also \cite{shiwang}) and, by properties of quasilinear elliptic operators, we have that $kF_{(u,v)}(\mu,u,v)+(1-k)F_{(u,v)}(\mu,\lambda,0)$ for $k\in(0,1)$  is Fredholm (see Sec 3. in \cite{shiwang}). Then, by Theorem 4.4 in \cite{shiwang}  (see also \cite{lopezgomez}, Theorem 6.4.3, page 188) the set $\mathcal{C}$ satisfies one of the following alternatives:
\begin{enumerate}[(i)]
    \item it is unbounded in $\mathbb{R}\times X_{\Omega}\times X_{\Omega_{1}}$, or
    \item it contains $(\mu_{*},\lambda,0)$, where $\mu_{*}\not=\mu_{\lambda}$ and $(\mu_{*},\lambda,0)$ also solves $F_{0}(\mu,u,v)=\Vec{0}$, or
    \item it contains a point $(\mu,\lambda + z,z)$, where $z\neq 0$ and $z$ is in the complement of $ker(F_{(u,v)}(\mu,\lambda,0))$ in $X_{\Omega}\times X_{\Omega_{1}}$.
    
\end{enumerate}
If we assume that (ii) holds and that $\mu<c/m$ then, by the same argument used to obtain (\ref{branchsolution}), we cannot have $\mu_{*}>\mu_{\lambda}$ since that would imply negative eigenvalues for the negative Laplacian with Neumann boundary conditions (see Eqn. \ref{equ2-6}). On the other hand, $\mu_{*}<\mu_{\lambda}$ is not possible since $\mu_{\lambda}$ is the smallest value for which positive solutions bifurcate. Now, let us assume that (iii) holds. Then, there is a point $(\mu,\lambda+z,z)\in\mathcal{C}$ with $z\neq 0$ and $\int_{\Omega_{1}}z\alpha_{\mu_{\lambda}}dx=0$. By the Sobolev's embedding theorem we can find a sequence $\{(\mu_{i},u_{i},v_{i})\}_{i={1}}^{\infty}$ in $\mathbb{R}\times C^{1}(\overline{\Omega})\times C^{1}(\overline{\Omega}_{1})$ with $u_{i}>0$ in $\Omega$, $v_{i}>0$ in $\Omega_{1}$ for all $i\in\mathbb{N}$ with $F_{0}(\mu_{i},u_{i},v_{i})=\Vec{0}$, and such that $\displaystyle\lim_{i\rightarrow\infty}(\mu_{i},u_{i},v_{i})=(\mu,\lambda+z,z)$. Then, we have that $(\lambda + z, z)$ is a non-negative solution of (\ref{system2}) with $\mu=\mu_{\lambda}$. By Lemma 2.2 in \cite{He2017} (see also, \cite{louni} and \cite{lintakagi}) and the Maximum Principle for quasilinear operators (see \cite{pucciserrin}) we must have
$$
\text{(1) }\lambda + z >0, z\equiv 0 \text{ or (2) }\lambda+z\equiv0, z\equiv0 \text{ or (3) } \lambda+z\equiv0, z>0.
$$
By assumption, neither $(1)$ nor $(2)$ hold. For the last case, the Maximum Principle implies that $\alpha_{\mu_{\lambda}}>0$ and therefore $\int_{\Omega_{1}}z\alpha_{\mu_{\lambda}}dx>0$, thus obtaining a contradiction. We have thus proved the following Theorem.
%

\begin{thm}\label{teorema1}Let $\lambda>0$. Then, there are positive solutions to the system (\ref{system2}) bifurcating from $\{(\mu,u,v)=(\mu,\lambda,0):\mu>0\}$ if and only if $0<\mu<\mu_{\lambda}$, forming a smooth curve given by  
\begin{align}
\{(\mu,u,v)=(\mu_{\lambda}(s),\lambda-s\alpha_{\mu_{\lambda}}(x)+o(|s|),s+o(|s|)):s\in(0,a)\}
 \end{align}
for some $a>0$ and such that $\mu_{\lambda}(0)=c\lambda/(1+m\lambda)$, $u(0)=\lambda$, $v(0)=0$. Furthermore, if $\mu>\mu_{\lambda}$ the system (\ref{system2}) has not positive solutions.
\end{thm}

\subsection{The linear diffusion case}
Now let us consider the system
\begin{equation}\label{system5}
\begin{aligned}
    \Delta u +\lambda u - u^{2} - \frac{b(x)uv}{1+mu} &= 0 \qquad\text{in}\quad \Omega \\ 
    \Delta v - \mu v + \frac{cuv}{1+mu} &= 0\qquad\text{in}\quad \Omega_{1}\\ 
    \partial_{n}u &= 0\qquad \text{on}\quad \partial{\Omega}, \\ \partial_{n}v &= 0 \qquad\text{on}\quad \partial\Omega_{1},\\ 
\end{aligned}
\end{equation}
which is identical to (\ref{system2}) but has the linear diffusion in $u$.
From the maximum principle stated in \cite{lou} (Proposition 2.2), any non-negative solutions to (\ref{system5}) are either positive, $(0,0)$, or $(\lambda,0)$. By letting $w=\lambda-u$, we define the operator $T:\mathbb{R} \times X_{\Omega}\times X_{\Omega_{1}}\rightarrow Y_{\Omega} \times Y_{\Omega_{1}}$ as 

\begin{equation}
T(\mu,w,v)=\begin{pmatrix}
\Delta w - \lambda w + w^{2} + \frac{b(x)(\lambda - w)v}{1+m(\lambda-w)}\\ \\
\Delta v - \mu v + \frac{c(\lambda - w)v}{1+m(\lambda - w)}
\end{pmatrix}^{T}.
\end{equation}
As before, we obtain the corresponding expressions for $T_{(w,v)}(\mu,w,v )$, $T_{(w,v)(w,v)}(\mu,w,v)$ and $T_{\lambda(w,v)}(\mu,w,v)$,
\begin{equation}
T_{(w,v)}(\mu,w,v)[\alpha,\beta]=
    \begin{pmatrix}
    \Delta\alpha - \lambda\alpha + 2w\alpha - \frac{b(x)\alpha v}{(1+m(\lambda - w))^{2}} + \frac{b(x)(\lambda - w)\beta}{1+m(\lambda - w)} \\ \\
    \Delta\beta - \mu\beta -\frac{c\alpha v}{(1+m(\lambda - w))^{2}} + \frac{c(\lambda - w)\beta}{1+m(\lambda - w)}
    \end{pmatrix}^{T}
\end{equation}{}
\begin{equation}
    T_{\mu(w,v)}(\mu,w,v)[\alpha,\beta] = (0,-\beta)
\end{equation}{}
\begin{equation}
    T_{(w,v)(w,v)}(\mu,w,v)[\alpha,\beta]^{2}=
    \begin{pmatrix}
    2\alpha^{2} - \frac{2b(x)}{(1+m(\lambda - w))^{2}}\alpha\beta- \frac{2mb(x)}{(1+m(\lambda - w))^{3}}\alpha^{2}v \\ \\
    -\frac{2c}{(1+m(\lambda - w))^{2}}\alpha\beta - \frac{2mc}{(1+m(\lambda - w))^{3}}\alpha^{2}v
    \end{pmatrix}^{T}.
\end{equation}{}
Therefore, by making $w=v=0$, 
\begin{equation}
T_{(w,v)}(\mu,0,0)[\alpha,\beta]=\begin{pmatrix}
\Delta\alpha - \lambda\alpha + b(x)\frac{\lambda}{1+m\lambda}\beta \\ \\
\Delta\beta - \mu\beta + c\frac{\lambda}{1+\lambda}\beta
\end{pmatrix}^{T}.
\end{equation}
By the same arguments used for (\ref{operator3}) we conclude that  $\mu_{\lambda}=c\lambda/(1+m\lambda)$ is the unique bifurcation along the curve $\Gamma_{u}$ from which positive solutions of (\ref{system5}) may  emerge. Notice also that
\begin{equation}
    T_{(w,v)(w,v)}(\mu_{\lambda},0,0)[\widetilde{\alpha}_{\mu_{\lambda}},\beta]^{2}=
    \begin{pmatrix}
    2\widetilde{\alpha}^{2}_{\mu_{\lambda}}-\frac{2b(x)}{(1+m\lambda)^{2}}\widetilde{\alpha}_{\mu_{\lambda}}\beta \\ \\
    -\frac{2c}{(1+m\lambda)^{2}}\widetilde{\alpha}_{\mu_{\lambda}}\beta
    \end{pmatrix}^{T},
\end{equation}
where
\[\widetilde{\alpha}_{\mu_{\lambda}}=(-\Delta + \lambda I)^{-1}\left[b(x)\lambda/(1+m \lambda)\right].\] 
Similarly as in the case of (\ref{derivativeeval}), but this time using $T_{(w,v)(w,v)}(\mu_{\lambda},0,0)\notin Range(T_{(w,v)}(\mu_{\lambda},0,0))$, we use formula I.6.3 in \cite{Hansjorg} to obtain
\begin{align}\label{derivativeeval2}
    \mu^{\prime}_{\lambda}(0) &= -\frac{1}{2}\frac{\langle T_{(w,v)(w,v)}(\mu_{\lambda},0,0)[\widetilde{\alpha}_{\mu_{\lambda}},\beta]^{2}, \pi_{1}\rangle}
{\langle T_{\mu(w,v)}(\mu_{\lambda},0,0)[\widetilde{\alpha}_{\mu_{\lambda}},\beta],\pi_{1} \rangle}\\
&=-\frac{c}{|\Omega_{1}|(1+m\lambda)^{2}}\int_{\Omega_{1}}\widetilde{\alpha}_{\mu_{\lambda}}dx < 0.
\end{align}
By defining $F_{1}:\mathbb{R}\times X_{\Omega}\times X_{\Omega_{1}}\rightarrow Y_{\Omega}\times Y_{\Omega_{1}}$ as
\begin{equation}
    F_{1}(\mu,u,v)=
    \begin{pmatrix}{}
    u-\lambda\\ \\
    v
    \end{pmatrix}^{T}-
    \begin{pmatrix}{}
    (-\Delta + I)^{-1}(u-\lambda +\lambda u - u^{2} - b(x)uv/(1+mu) )\\
    \\
    (-\Delta + I)^{-1}_{\Omega_{1}}(v - \mu v + cuv/(1+mu))
    \end{pmatrix}^{T}
\end{equation}
we can use Shi and Wang's result, \cite{shiwang}, (see also \cite{lopezgomez}), to get an analogous statement to Theorem \ref{teorema1}. Summarizing, we have the following result.

\begin{thm}\label{teorema2}Let $\lambda>0$. Then, there are positive solutions to the system (\ref{system2}) bifurcating from $\{(\mu,u,v)=(\mu,\lambda,0):\mu>0\}$ if and only if $0<\mu<\mu_{\lambda}$, forming a smooth curve given by  
\begin{align}
\{(\mu,u,v)=(\mu_{\lambda}(s),\lambda-s\widetilde{\alpha}_{\mu_{\lambda}}(x)+o(|s|),s+o(|s|)):s\in(0,a)\}
 \end{align}
for some $a>0$ and such that $\mu_{\lambda}(0)=c\lambda/(1+m\lambda)$, $u(0)=\lambda$, $v(0)=0$. Furthermore, if $\mu>\mu_{\lambda}$ the system (\ref{system5}) has not positive solutions.
\end{thm}

\section{Numerical verification}
Bifurcation curves in the $v-\mu$ plane, for the system (\ref{system2}) and the linear counterpart, were computed for three values of $\lambda$ and are presented in Figure 1. The numerical results obtained are in agreement with the theoretical findings in the previous Section. The bifurcation points coincide for both systems and the emerging curves of positive solutions are virtually identical for values of $\mu$ that are below but close to the bifurcation point. However, it is observed that as the values of $\mu$ move toward 0, the curve associated to the nonlinear diffusion eventually starts to increase much faster than that of its linear counterpart.    

\section{Conclusions and discussion}
In this note we introduce a simple nonlinear diffusion mechanism in a prey population to model a plausible adaptation response that counteracts intraspecific competition for resources. The spatial domain for the model contains a refuge zone that excludes predators presence. Predator saturation on prey consumption is also considered and included via a Holling type II function. It is reasonable to imagine this scenario as a simplified approximation to biological control or conservation problems where the question ``\textit{how does density-dependent diffusion in the prey affect the dynamics of the system and compares to linear diffusion?}" is relevant to a modeler of such complex contexts. In this paper we provide a partial answer to this general question for a very simple case of nonlinear diffusion. Our theoretical arguments involve a novel adaptation of a maximum principle to the nonlinear case and make use of the classical results in bifurcation theory to show the existence of positive solutions at the steady state. The analysis is complemented by the numerical computation of the bifurcation curves for the nonlinear and linear diffusion cases.
In addition, it can be shown that the observed bifurcation is transcritical, the details of the proof can be found in \cite{LeoTranscritic}. 

Our study complements the literature on the theme, see for example \cite{chen,DuDu,LiDong,lou,oeda,wangh,wangjandwangy,He2017,zhangrongzhang}, and opens some new questions. For instance, it would be of use to find under which circumstances the theoretical framework can be extended to more general forms of density dependence, and if so, how do they compare with the linear diffusion case. Although interesting, these are out of the scope of this note.

\begin{figure}
\centering
\includegraphics[scale=0.22]{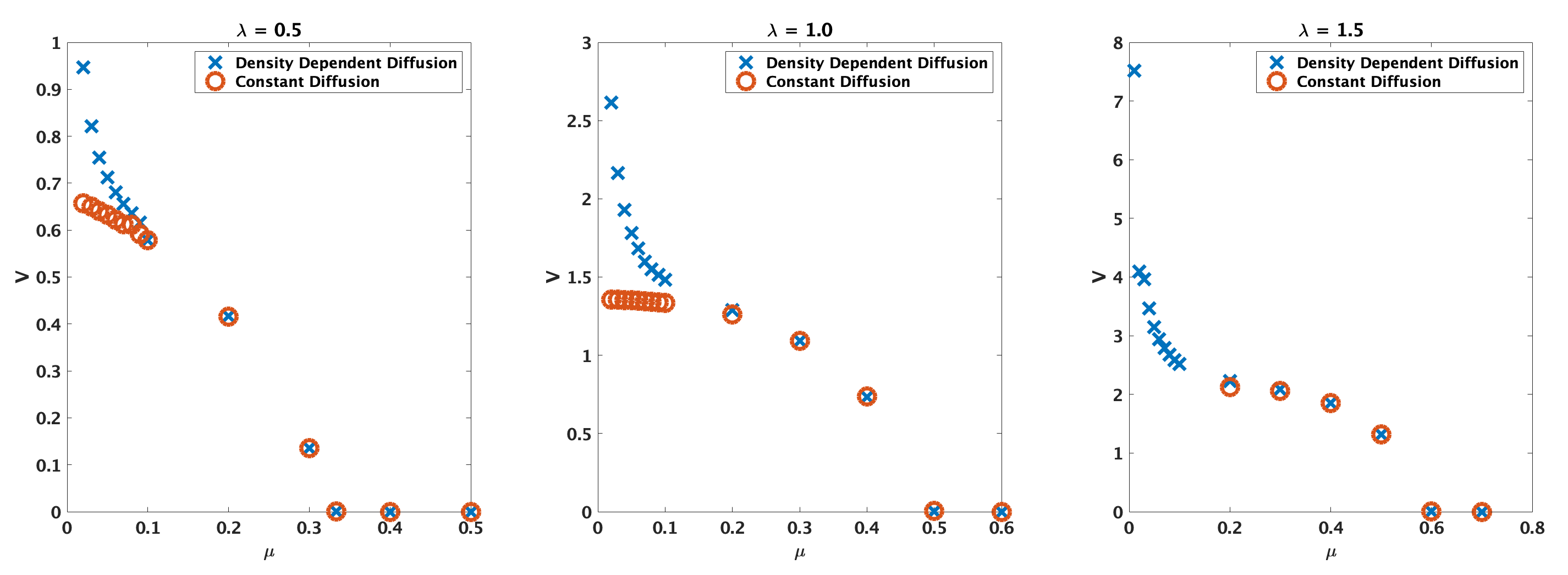}
\caption{Bifurcation curves for the density dependent and linear (constant) diffusion, with the predator population $v$ on the vertical axis and the parameter $\mu$ on the horizontal axis. The curves for the system (1.5) are in blue ($\times$), and for its linear counterpart in orange ($\circ$). From left to right, the values of $\lambda$ are 0.5, 1.0, and 1.5, with $c=m=1$ for the three panels. This gives the corresponding bifurcation points at $\mu=$1/3, 1/2, and 3/5. Close to the bifurcation points the curves of positive solutions for the nonlinear and linear diffusion are virtually the same but, as the value of $\mu$ moves toward 0, the former increases faster.}
\label{fig:method}
\end{figure}

\textbf{Acknowledgements.} The authors are grateful to an anonymous reviewer for pointing out a gap in the proof of Theorem 1, and for suggesting the references \cite{shiwang} and \cite{lopezgomez}. Many thanks to Paul Sacks and Z.Q. Wang for useful comments. LRQ is grateful to the Department of Mathematics and Statistics for support.


\newpage
\bibliography{REFS}
}\end{document}